\documentclass[11pt]{amsart}

\usepackage{amsfonts}
\usepackage{amsmath}
\usepackage{amssymb, esint}
\usepackage{amsthm,color}
\usepackage{enumerate}
\usepackage{mathtools}
\usepackage{cite}

\newcommand{\R}{{\mathbb R}}

\newcommand{\N}{{\mathbb N}}

\def\0{{\mathbf 0}}

\newcommand{\e}{\varepsilon}
\newcommand{\vp}{\varphi}

\newcommand{\supp}{\operatorname{spt}}

\def\mean#1{\mathchoice%
          {\mathop{\kern 0.2em\vrule width 0.6em height 0.69678ex depth -0.58065ex
                  \kern -0.8em \intop}\nolimits_{\kern -0.4em#1}}%
          {\mathop{\kern 0.1em\vrule width 0.5em height 0.69678ex depth -0.60387ex
                  \kern -0.6em \intop}\nolimits_{#1}}%
          {\mathop{\kern 0.1em\vrule width 0.5em height 0.69678ex
              depth -0.60387ex
                  \kern -0.6em \intop}\nolimits_{#1}}%
          {\mathop{\kern 0.1em\vrule width 0.5em height 0.69678ex depth -0.60387ex
                  \kern -0.6em \intop}\nolimits_{#1}}}

\theoremstyle{plain}
\newtheorem{thm}{Theorem}[section]

\newtheorem{lem}[thm]{Lemma}
\newtheorem{prop}[thm]{Proposition}

\theoremstyle{defn}
\newtheorem{defn}[thm]{Definition}

\theoremstyle{rem}

\newtheorem{rem}[thm]{Remark}

\numberwithin{equation}{section}

\title[]{ Almost minimizers to a   transmission problem  for $(p,q)$-Laplacian}

\author{Sunghan Kim}
\address{Department of Mathematics, Uppsala University, S-751 06 Uppsala, Sweden}
\email{sunghan.kim@math.uu.se}

\author{Henrik Shahgholian}
\address{Department of Mathematics, Royal Institute of Technology, 100 44 Stockholm, Sweden}
\email{henriksh@kth.se}

\thanks{ H. Shahgholian was supported in part by Swedish Research Council. This project was finalized during the program Geometric aspects of nonlinear PDE at Institute Mittag Leffler, Stockholm.}

\begin{document}

\maketitle

\begin{abstract}
This paper concerns almost minimizers of the functional 
$$
J(v,\Omega) = \int_\Omega \left( |D v^+|^p +  |D v^-|^q \right) dx,
$$
where $1<p \neq q< \infty$ and $\Omega$ is a bounded domain of $\R^n$, $n\geq 1$. We prove the universal H\"older regularity of local $(1+\e)$-minimizers, when $\e$ is universally small. Moreover, we prove almost Lipschitz regularity of the local $(1+\e)$-minimizers, when $|p-q|\ll 1$ and $\e\ll 1$. 

\end{abstract}

\tableofcontents


\section{Introduction}\label{section:intro}

In this paper, we study regularity properties of  almost minimizers to the functional 
\begin{equation}\label{eq:J}
J(u,\Omega) \equiv J_{p,q}(u,\Omega) := \int_\Omega (|Du^+|^p + |Du^-|^q)\,dx,
\end{equation}
where $\Omega\subset\R^n$ is a bounded domain and $1<p , q <\infty$. Our primary goal is to prove a universal H\"older estimate for the almost minimizers. We shall also study various scenarios, on the relation between $p$ and $q$, to see if the regularity can be improved. In particular, we aim at proving almost Lipschitz regularity provided that $p$ and $q$ are close to each other.

The notion of local $K$-minimizers is given as follows. 

\begin{defn}[Local $K$-minimizers]
Let $K\geq 1$ be a constant. We shall call $u\in W_{loc}^{1,p\wedge q}(\Omega)$ a local $K$-minimizer of the functional $J$, if for any cube $Q\subset\Omega$, $J(u, Q) < \infty$, and 
\begin{equation}\label{eq:defn}
J(u,Q) \leq K J(v,Q),
\end{equation}
for any $v\in u + W_0^{1,p\wedge q}(Q)$ such that $J(v,Q) < \infty$. 
\end{defn}
In the course of this paper, we shall be interested in the case $K = 1 + \e$, for some small $\e>0$. We remark that our analysis does not change, as one replaces cubes with balls in the above definition. However, it is worth mentioning that the notion with cubes is in general not equivalent to that with balls, unless $K = 1$, and local $K$-minimizers with cubes are known to be less restrictive; see \cite[Example 6.5]{Giu}. 

In the framework of standard functionals (i.e., those without break across some level set), the universal H\"older regularity is established for quasi-minimzers (those with $K>1$ any, and $Q$ in \eqref{eq:defn} replaced with $\supp (u-v)$), as the essential arguments for the proof of the H\"older regularity for exact minimizers remain unchanged upon the extension; see \cite{Giu}. In contrast, thanks to the particular break across the zero-level set in $J_{p,q}$, many important steps in the proof of \cite[Theorem 1.2]{CKS} for the H\"older regularity of exact minimizers to our functional $J_{p,q}$ are destroyed when applied to quasi-minimzers. Still, we were able to extend the argument to $(1+\e)$-minimizers, when $\e$ is universally small. 

\begin{thm}\label{thm:Ca}
There are constants $\e>0$ and $\sigma\in(0,1)$, depending only on $n$, $p_+$, and $p_-$, such that if $u\in W^{1,p_+\wedge p_-}(Q_2)$ is a local $(1+\e)$-minimizer of $J_{p_+,p_-}$, then $u^\pm \in C_{loc}^{0,\sigma_\pm}(Q_1)$ with $\sigma_+ = \sigma$, $\sigma_- = 1- (1-\sigma)\frac{p_-}{p_+}$, and 
$$
[ u^\pm ]_{C^{0,\sigma_\pm}(Q_1)} \leq c\left[\int_{Q_2} ((u^+)^{p_+} + (u^-)^{p_-})\, dx\right]^{\frac{1}{p_\pm}}, 
$$
where $c$ depends only on $n$, $p_+$, and $p_-$. 
\end{thm}

We remark that the above theorem also shows the exact relation between the H\"older exponents for each phase; this was not contained in the authors earlier collaboration \cite[Theorem 1.2]{CKS} with M. Colombo. Our proof involves a careful extension of the main ingredients for \cite[Theorem 1.2]{CKS} to local $(1+\e)$-minimizers, and a compactness argument. 

A key feature of local $(1+\e)$-minimizers, $\e\geq 0$, for the functional $J_{p,q}$ is that the positive and negative phase scales differently from each other. Namely if $u$ is a local $(1+\e)$-minimizer in $Q_2$, then one needs $\| u^+ \|_X$ comparable with $\|u^- \|_X^{q/p}$, with $X = L^p(Q_1)$ or $L^\infty(Q_1)$. As for the case of the local minimizers, i.e., $\e = 0$, the comparability was proved by a Harnack inequality argument \cite[Lemma 3.8, Corollary 3.9]{CKS}, which played an essential role in the proof of their universal H\"older regularity \cite[Theorem 1.2]{CKS}. 

The main difference, which also amounts to the challenges here, for the case of local $(1+\e)$-minimizers, $\e>0$, is the lack of such a Harnack inequality argument. More fundamentally, local $(1+\e)$-minimizers do not possess the subsolution properties as opposed to local minimizers (see \cite[Lemma 3.4]{CKS}). One of the consequences is that the basic estimates for one phase, such as the Cacciopoli inequality (Lemma \ref{lem:caccio}) and the comparison lemma (Lemma \ref{lem:comp}) for local $(1+\e)$-minimizers, involve an additional $\e$-factor of the other phase. Hence, our main task here is to effectively control the additional $\e$-term, which amounts to some technical difficulties. It is worthwhile to mention that  the absence of the Harnack inequality argument is overcome by a careful compactness argument, by which both phases, although scaled differently, survive at the limit. The latter part is new, to the best of the authors' knowledge, and can be applied to a wider range of problems. 

Our second result is about the almost Lipschitz regularity for local $(1+\e)$-minimizers for the functional $J_{p,q}$, when $|p-q| \ll 1$ and $\e \ll 1$. 

\begin{thm}\label{thm:C01-} 
Let $1 < p_+  <\infty$ and $\sigma \in (0,1)$ be given. Then there exist $\e,\delta > 0$, depending only on $n$, $p_+$ and $\sigma$, such that for any $p_-\in(p_+-\delta,p_+ + \delta)$ and any local $(1+\e)$-minimizer $u\in W^{1,p_+\wedge p_-}(Q_2)$ of $J_{p_+,p_-}$, one has $u^\pm\in C^{0,\sigma_\pm}(Q_1)$, with $\sigma_+ = \sigma$, $\sigma_- = 1 -  (1-\sigma)\frac{p_-}{p_+}$, and 
$$
[u]_{C^{0,\sigma_\pm}(Q_1)} \leq c \left[ \int_{Q_2} ((u^+)^{p_\pm} + (u^-)^{p_-})\,dx\right]^{\frac{1}{p_\pm}},
$$
where $c$ depends only on $n$, $p_+$ and $\sigma$. 
\end{thm}

A similar statement is proved in \cite{AT} for uniformly elliptic functionals when governing conductivity matrices are close with each other; \cite{AT} however considers local minimizers (i.e., $\e = 0$) only. Our problem is philosophically the same, as the limit case is clean, thus possess better regularity. On the technical level, 
our argument needs slightly more care
 than that of \cite[Theorem 7.1]{AT}, as the proof for the growth of the functional $J_{p,q}$ changes as $(p,q)$ varies. Moreover, one needs to make sure that the argument works well regardless of the relation between $p$ (or $q$) and the dimension $n$. These are all rigorously treated in Sect.~\ref{sec:C01-}. 

Recently, free boundaries for almost minimizers are investigated in various settings, see e.g., \cite{DET}, \cite{DS}, and \cite{DJS} to mention a few. There is a possibility of extending the approach with viscosity solutions employed in \cite{DS}, but it is beyond the scope of this paper. It would be already interesting to extend the result for the clean case, $p = q$. 

In \cite{CKS}, the authors analyze the free boundary of local minimizers for $J_{p,q}$, using the measure $\Delta_p u^+$, which is nonnegative and supported on the free boundary, $\partial\{u > 0\}$(=$\partial\{u < 0\}$). This is mainly due to the subsolution property of $u^+$, which is no longer valid for almost minimizers. The same issue appears in the case of the two-phase Alt-Caffarelli functional (see \cite[Section 4]{DET}), which is resolved by the NTA property of the free boundary and a clever use of barriers. The NTA property was obtained there by the use of the ACF monotonicity formula, which is absent in our regime. The construction of the barriers and the comparison with the almost minimizers require some regularity of the free boundary, which in the case of \cite{DET} was the NTA property. However, in our problem, none of these seems to be analogously carried out. For this reason, we leave out the analysis of the free boundary for our almost minimizers to the interested reader. 

The paper is organized as follows. In Section \ref{sec:tech}, we collect some technical tools to prepare the proof of Theorem \ref{thm:Ca}. In Section \ref{sec:Ca}, we prove Theorem \ref{thm:Ca}. In Section \ref{sec:C01-}, we prove Theorem \ref{thm:C01-}. 

We follow the standard notation and terminology. In particular, $n$ denotes the dimension of the underlying space, and there is no restriction other than $n\geq 1$. By $Q_r(x_0)$, we denote the cube centered at $x_0$ with side-length $r$, i.e., $Q_r(x_0) := \{ x \in \R^n: |x_i - x_{0i}| < r,\, 1\leq i\leq n\}$. For simplicity, we set $Q_r := Q_r(0)$. Given a set $A\subset\R^n$, by $|A|$ we denote the Lebesgue measure of $A$. The function spaces $C^{0,\sigma}$ and $W^{1,p}$ are standard H\"older and Sobolev spaces, and $C_{loc}^{0,\sigma}$, $W_{loc}^{1,p}$ are their local versions.


\section{Technical Tools}\label{sec:tech}

In this section, we shall present and verify some technical tools, most of which generalize those appeared in \cite[Sect. 3--4]{CKS}. The main goal of this section is to prove the following proposition, which roughly tells us that negative values cannot penetrate the interior if a local $(1+\e)$-minimizer attains large positive values in most of the domain. Let us remark that this proposition corresponds to \cite[Proposition 4.2]{CKS} for the case of minimizers. 
The main difference here is that $(1+\e)$-minimizers do not possess in general the subsolution properties. Here we exploit the techniques to circumvent this issue. Unless stated otherwise, the constant $c$ throughout this section is a positive constant that may differ at each occurrence, and will depend at most on $n$, $p$, and $q$. Moreover, the parameter $\e$ will be a small constant, whose smallness is determined solely by $n$, $p$, and $q$. 

\begin{prop}\label{prop:pos}
There exist $\e>0$ and $\mu>0$, depending only on $n$, $p$, and $q$, such that if $u\in W^{1,p\wedge q}(Q_1)$ is a local $(1+\e)$-minimizer of the functional $J$, satisfying
$$ 
\int_{Q_1} ((u^+)^p + (u^-)^q)\,dx \leq 1,\qquad |\{ u\leq 1/2\}\cap Q_1|\leq \e, 
$$
then $u > 0$ a.e.\ in $Q_\mu$. 
\end{prop}

The proof for this proposition will be postponed to the end of this section. Let us begin with the Cacciopoli-type inequality.

\begin{lem}\label{lem:caccio}
Let $u\in W^{1,p\wedge q}(Q_2)$ be a local $(1+\e)$-minimizer of the functional $J$. There exists $\bar\e \in (0,1)$, depending only on $n$, $p$, and $q$, such that if $\e \leq \bar\e$, then 
\begin{equation}\label{eq:caccio}
\int_{Q_1} |Du^+|^p \,dx \leq c \int_{Q_2} ((u^+)^p + \e (u^-)^q)\,dx,
\end{equation}
where $c$ depends only on $n$, $p$, and $q$. 
\end{lem}

\begin{proof} 
Fix $r,R$ with $1<r<R<2$, and choose any $s,t$ with $r<s<t<R$. Let $\eta \in C_c^1(Q_t)$ be a cutoff function such that $\eta\equiv 1$ in $Q_s$, $|D\eta|\leq \frac{2c}{t-s}$ in $Q_t$, and $\supp(\eta)\subset Q_{(t+s)/2}$. Set $w := (1 - \eta)u^+ - u^- \in W^{1,p\wedge q}(Q_t)$. Since $w^+ = (1-\eta)u^+$, $w^- = u^-$, and $\supp(u-w) \subset \supp(\eta) \subset Q_{(t+s)/2}$, we derive from the $(1+\e)$-minimizerslity of $u$ for $J_{p,q}$ in $Q_t$ that  
$$
\int_{Q_r} |Du^+|^p \,dx \leq (1+\e) \int_{Q_t} |D((1-\eta)u^+)|^p\,dx + \e \int_{Q_t} |Du^-|^q\,dx. 
$$
Applying H\"older's inequality and Young's inequality, and then using $\supp(\eta)\subset Q_{(t+s)/2}$ and $|D\eta|\leq c/(t-s)$, we deduce that 
$$
\int_{Q_s} |Du^+|^p\,dx \leq c\int_{Q_t} \left(\frac{(u^+)^p}{(t-s)^p} + \e |Du^-|^q\right) dx+ c\e \int_{Q_t} |Du^+|^p \,dx.
$$
Since this part is by now standard, we omit the details. Note that the last display holds for all $s,t$, $r<s<t<R$. Hence, choosing $\e$ small enough such that $c\e < \frac{1}{2}$, we can employ the standard iteration lemma \cite[Lemma 6.1]{Giu} to derive that 
\begin{equation}\label{eq:caccio-re}
\int_{Q_r} |Du^+|^p\,dx \leq c\int_{Q_R} \left(\frac{(u^+)^p}{(R-r)^p} + \e |Du^-|^q\right) dx. 
\end{equation}

Now replace $Q_R$ in the right-hand side with $Q_{(R+r)/2}$, and then apply the same argument above to $(-u)$ with $Q_r$ replaced with $Q_{(R+r)/2}$; note that $(-u)$ is a local $(1+\e)$-minimizer of $J_{q,p}$ in place of $J_{p,q}$. Then we may proceed as follows, 
$$
\begin{aligned}
\int_{Q_r} |Du^+|^p\,dx &\leq c\int_{Q_{(R+r)/2}} \left(\frac{(u^+)^p}{(R-r)^p} + \e |Du^-|^q\right) dx \\ 
&\leq c\int_{Q_R}  \left(\frac{(u^+)^p}{(R-r)^p} + c\e \frac{(u^-)^q}{(R-r)^q}\right)dx   + c^2\e^2\int_{Q_R}  |Du^+|^p\, dx. 
\end{aligned}
$$
Recall that $r,R$ were any numbers between $1$ and $2$. Hence, taking $\e$ smaller if necessary such that $c^2\e^2 < \frac{1}{2}$, we can make use of the iteration lemma once again to arrive at \eqref{eq:caccio}. 
\end{proof}

\begin{rem}
In what follows, we shall always assume that $\e < \bar\e$, with $\bar\e$ as in Lemma \ref{lem:caccio}. 
\end{rem} 

Let us remark that the above Cacciopoli inequality is too weak to bring forth a local $L^\infty$-estimate. Besides, local quasi-minimizers are not necessarily bounded, even for functionals under standard growth condition (of course, only if $p \leq n$). Nevertheless, with the aid of the Cacciopoli inequality above, we shall observe that the blowup rate of local $(1+\e)$-minimizers can be made arbitrarily small, for small $\e$, in case $p \leq n$. 

\begin{lem}\label{lem:blow}
Let $u\in W^{1,p\wedge q}(Q_1)$ be a local $(1+\e)$-minimizer of the functional $J$. Suppose that 
\begin{equation}\label{eq:blow-asmp}
\| u^+ \|_{L^p(Q_1)} \leq 1,\quad \sup_{r\in(0,1)} \frac{\| u^-\|_{L^q(Q_r)}}{ r^{1-\frac{p}{q}}\| u^+\|_{L^p(Q_r)}^{\frac{p}{q}}} \leq \kappa,
\end{equation}
for some constant $\kappa>0$. Then for any $\delta>0$, there exists a positive constant $\e_{\kappa,\delta}$, depending only on $n$, $p$, $q$, $\kappa$ and $\delta$, such that if $\e \leq \e_{\kappa,\delta}$, then 
$$
\sup_{r\in(0,1)} \frac{1}{r^{n - \delta p}} \int_{Q_r} (u^+)^p\,dx \leq c_{\kappa,\delta},
$$
where $c_{\kappa,\delta}$ depends only on $n$, $p$, $q$, $\Lambda$, $\delta$ and $\kappa$.
\end{lem}

\begin{proof}
We remark that the conclusion is trivial for $p> n$, due to the Sobolev embedding theorem. Henceforth, we shall assume that $1 < p\leq n$. 

Let  $\kappa$ and  $\delta$ be arbitrary positive constants, and 
 suppose  the conclusion of the lemma is false. Then for each $j=1,2,\cdots$, one can find some positive constant 
 $\e_j \searrow 0$, and a local $(1+\e_j)$-minimizer $u_j \in W^{1,p\wedge q}(Q_1)$ of the functional $J$, such that 
$$
\| u_j^+ \|_{L^p(Q_1)} \leq 1,\quad \sup_{r\in(0,1)} \frac{\| u_j^-\|_{L^q(Q_r)}}{ r^{1-\frac{p}{q}}\| u_j^+\|_{L^p(Q_r)}^{\frac{p}{q}}} \leq \kappa,
$$
but
$$
S_j = \sup_{r_j \leq r\leq 1} \frac{1}{r^{n - \delta p}} \int_{Q_r} (u_j^+)^p\,dx \to \infty,
$$
for some constant $r_j\in(0,1)$. In order to have $S_j\to \infty$ to be compatible with $\|u_j^+\|_{L^p(Q_1)} = 1$, we must have $r_j \to 0$. 

Consider an auxiliary function $v_j : Q_{r_j^{-1}} \to \R$, defined by 
$$
v_j(y) = \frac{u_j^+(r_jy)}{r_j^{-\frac{n}{p}}\| u_j^+\|_{L^p(Q_{r_j})}} - \frac{u_j^-(r_jy)}{r_j^{1-\frac{p}{q} - \frac{n}{q}} \| u_j^+\|_{L^p(Q_{r_j})}^{\frac{p}{q}}}. 
$$
One  easily verifies  that $v_j\in W^{1,p\wedge q}(Q_{r_j^{-1}})$ is a local $(1+\e_j)$-minimizer of the functional $J$, and  
\begin{equation}\label{eq:vj+-Lp}
\sup_{1\leq R\leq r_j^{-1}} \frac{1}{R^{n - \delta p}} \int_{Q_R} (v_j^+)^p\,dy = 1, 
\end{equation}
where the supremum is in particular attained at  $R= 1$. Furthermore, from the assumption between the growth of $(u_j^-)^q$ and $(u_j^+)^p$ (with the correct scaling shown in the statement above), we obtain from \eqref{eq:blow-asmp} and \eqref{eq:vj+-Lp} that 
$$
\begin{aligned}
    \int_{Q_R} (v_j^-)^q\,dy & = \frac{1}{r_j^{n + q - p - n} \| u_j^+ \|_{L^p(Q_{r_j})}^p} \int_{Q_{Rr_j}} (u_j^-)^q\,dx  \\
    & \leq \frac{\kappa^q (Rr_j)^{q-p} \| u_j^+ \|_{L^p(Q_{Rr_j})}^p}{r_j^{q - p} \| u_j^+ \|_{L^p(Q_{r_j})}^p} \\
    & = \kappa^q R^{q-p} \int_{Q_R} (v_j^+)^p\,dy \\
    & \leq \kappa^q R^{n+q - (1+\delta)p},
\end{aligned}
$$
which holds for all $R \geq 1$. Therefore, we have  
\begin{equation}\label{eq:vj--Lq}
\begin{aligned}
\sup_{1\leq R\leq r_j^{-1}} \frac{1}{R^{n + q - (1+\delta)p}} \int_{Q_R} (v_j^-)^q\,dy \leq \kappa^q. 
\end{aligned}
\end{equation}

Due to Lemma \ref{lem:caccio}, along with \eqref{eq:vj+-Lp} and \eqref{eq:vj--Lq}, 
\begin{equation}\label{eq:Dvj-Lpq}
\int_{Q_R}( |Dv_j^+|^p + |Dv_j^-|^q )\,dx \leq c R^{n - (1 + \delta)p},
\end{equation}
where $c$ depends only on $n$, $p$ and $q$, whenever $2Rr_j \leq 1$. By the Sobolev embedding theory, there exists a function $v\in W_{loc}^{1,p\wedge q}(\R^n)$ with $v^+ \in W_{loc}^{1,p}(\R^n)$ and $v^- \in W_{loc}^{1,q}(\R^n)$ such that  $v_j^+ \to v^+$ and $v_j^- \to v_j^-$ weakly in $W_{loc}^{1,p}(\R^n)$ and respectively $W_{loc}^{1,q}(\R^n)$, after extracting a subsequence if necessary; we shall denote this subsequence by $v_j$, for brevity. 
The lower semicontinuity of the functional (see \cite[Proposition 2.2]{CKS}) along with the weak convergence
implies that $v\in W^{1,p\wedge q}(B_R)$ is a minimizer of the functional $J$. Since $v_j^+ \to v^+$ strongly in $L^p(B_R)$ and $v_j^-\to v^-$ strongly in $L^q(B_R)$, letting $j\to\infty$ in \eqref{eq:vj+-Lp} yields that
\begin{equation}\label{eq:v0+-Lp}
\sup_{R\geq 1} \frac{1}{R^{n-\delta p}} \int_{Q_R} (v^+)^p\,dy = 1.
\end{equation}
However, since $v$ is a minimizer of the functional $J$, by \cite[Lemma 3.4]{CKS}, $v^+$ is a weak $p$-subsolution. As a result, the local $L^\infty$-estimates \cite[Theorem 7.3]{Giu} applies to $v^+$, which along with \eqref{eq:v0+-Lp} yields 
$$
\| v^+\|_{L^\infty(Q_R)} \leq \frac{c}{R^\delta}. 
$$
Hence, letting $R\to\infty$ yields that $v^+ = 0$ a.e.\ in $\R^n$. This yields a contradiction against \eqref{eq:v0+-Lp}. 
\end{proof} 

We also have a growth estimate for the $p$-th Dirichlet energy of the positive phase. The idea is the same as in \cite[Lemma 3.6]{CKS}, which is based on some approximation by positive $p$-harmonic functions of the positive phase of local quasi-minimizers, in terms of the size of the negative phase. 

\begin{lem}\label{lem:growth}
Let $u\in W_{loc}^{1,p\wedge q}(Q_2)$ be a local $(1+\e)$-minimizer of the functional $J$, and $v \in u^+ + W_0^{1,p}(Q_1)$ be the $p$-harmonic function. Then 
$$
0\leq \int_{Q_1} (|Du^+|^p - |Dv|^p)\,dx \leq c\int_{Q_2}((u^-)^q + \e |Du^+|^p)\,dx,
$$
and
$$
\int_{Q_r} |Du^+|^p\,dx \leq c \int_{Q_1} ((r^n + \e) |Du^+|^p + (u^-)^q)\,dx,\quad\forall r\in(0,1),
$$
where $c$ depends only on $n$, $p$ and $q$. 
\end{lem}

\begin{proof}
The proof is essentially the same as that of \cite[Lemma 3.6]{CKS}. The additional term $\e \int_{Q_2} |Du^+|^p\,dx$ appears due to the different Cacciopoli inequality; more exactly, we use \eqref{eq:caccio-re} with $u$ replaced with $-u$. We shall not repeat this argument here. 
\end{proof}

The following lemma corresponds to \cite[Lemma 3.7]{CKS}. The key ingredient of the proof there is the Poincar\'e inequality, and Lemma \ref{lem:growth}, which corresponds to \cite[Lemma 3.6]{CKS}. As noted above, Lemma \ref{lem:growth} differs from \cite[Lemma 3.6]{CKS} by the additional term, $\e\int_{Q_2} |Du^+|^p\,dx$. However, this does not make any difference in the proof of the lemma below. Thus, we shall skip the proof.

\begin{lem}[Essentially due  to {\cite[Lemma 3.7]{CKS}}] 
\label{lem:msr} 
Let $u\in W^{1,p\wedge q}(Q_4)$ be a local $2$-minimizer for the functional $J$, satisfying 
$$
\mean{Q_4} (u^+)^p \,dx = 1,\quad \mean{Q_4} ( (u^-)^q +  |Du^+|^p)\,dx  \leq \e ,
$$
for some $\e>0$. Then
$$
|\{ u\leq 1/2\}\cap Q_1|\leq c\e,
$$
where $c$ depends only on $n$, $p$ and $q$. 
\end{lem}

Let us prove Proposition \ref{prop:pos} with additional assumptions that $\| u^-\|_{L^q(Q_1)}$ and $\|Du^+\|_{L^p(Q_!)}$ are sufficiently small. The proof follows the idea of that of \cite[Lemma 4.3]{CKS}, with some modifications addressing the lack of subsolution properties of each phase. 

\begin{lem}\label{lem:pos}
There exists $\e>0$, depending only on $n$, $p$ and $q$, such that if $u\in W^{1,p\wedge q}(Q_4)$ is a local $(1+\e)$-minimizer of the functional $J$, satisfying
$$
\mean{Q_4} (u^+)^p\,dx = 1,\quad\mean{Q_4} ((u^-)^q + |Du^+|^p)\,dx \leq\e,
$$
then $u> 0$ a.e.\ in $Q_1$. 
\end{lem} 

\begin{proof}
Let us consider the case $q < n$ first. Following the proof of \cite[Lemma 4.3]{CKS}, we obtain that for $\sigma\in(0,1)$,
\begin{equation}\label{eq:u-Lpq-gr}
\mean{Q_r} \left(\frac{(u^-)^q}{r^q} + |Du^+|^p\right) dx \leq c\e r^{-(1-\sigma )p} \mean{Q_r} (u^+)^p\,dx,\quad\forall r\in(0,1),
\end{equation}
where $c$ depends only on $n$, $p$, $q$ and $\sigma$. The proof is essentially the same, as Lemma \ref{lem:growth} and \ref{lem:msr} replace  \cite[Lemma 3.5--3.7]{CKS}, which are the key ingredients of the proof there; moreover Lemma \ref{lem:caccio} replaces the usual Cacciopoli inequality for weak $q$-subsolutions. These lemmas have additional $\e$-term, which arise from the $(1+\e)$-local minimizerslity of $u$, but this does not contribute any major difference from the proof for \cite[Lemma 4.3]{CKS}. Hence, we shall omit the details. 

We observe that due to \eqref{eq:u-Lpq-gr} (as well as the assumption $\mean{Q_4} (u^+)^p\,dx = 1$), the hypothesis of Lemma \ref{lem:blow} is satisfied (with $\kappa =  1 > \e r^{\sigma p}$). Thus, choosing $\e\leq \e_\delta$ with $\e_\delta$ as in Lemma \ref{lem:blow} with $\delta < \sigma$, we deduce
\begin{equation}\label{eq:u+-Lp-gr}
\mean{Q_r} (u^+)^p \,dx \leq c r^{- \delta p},\quad\forall r\in(0,1). 
\end{equation}
Inserting \eqref{eq:u+-Lp-gr} into \eqref{eq:u-Lpq-gr} yields that 
\begin{equation}\label{eq:Du+-Lp-gr-re}
\mean{Q_r} |Du^+|^p\,dx \leq c\e r^{-(1 - (\sigma -\delta))p},\quad\forall r\in(0,1);
\end{equation} 
now $c$ depends only on $n$, $p$, $q$, $\sigma$ and $\delta$. Let us remark that this step does not appear for the case of minimizers \cite[Lemma 4.3]{CKS} because  for the latter case we can use the subsolution property \cite[Lemma 3.4]{CKS} for $u^+$ to obtain its local boundedness.

The growth estimate in \eqref{eq:Du+-Lp-gr-re} is obtained by choosing $\e$ 
 sufficiently small. Taking $\e$ even smaller if necessary, we may repeat the above argument around any point $z\in Q_1$, and obtain
$$
\mean{Q_r(z)} |Du^+|^p \,dx \leq c\e r^{-(1-(\sigma -\delta))p},\quad\forall r\in(0,1),\,\forall z\in Q_1,
$$
possibly with a larger constant $c$. Therefore, by Morrey's lemma, we deduce that $u^+ \in C^{0,\sigma-\delta}(Q_1)$ and
\begin{equation}\label{eq:u+-Ca-re}
[u^+]_{C^{0,\sigma-\delta}(Q_1)} \leq c\e^{\frac{1}{p}}.
\end{equation}
Finally, by Lemma \ref{lem:msr}, $|\{ u \leq \frac{1}{2}\}\cap Q_1| \leq c\e$. Hence, with $c\e \leq 2^{-2n -1}$, we have $|\{ u > \frac{1}{2}\}\cap Q_1| > 0$, which now implies via \eqref{eq:u+-Ca-re} that 
$$
\inf_{Q_1} u^+ \geq \frac{1}{2} - c \e^{\frac{1}{p}} > 0,
$$
provided that we choose $\e$ even smaller. Note that the smallness condition for $\e$ at this stage can be determined solely by $n$, $p$ and $q$, by for instance selecting $\sigma = \frac{1}{2}$ and $\delta = \frac{1}{4}$. This finishes the proof for the case $q < n$. 

The case for $q\geq n$ can be treated similarly, following the proof of \cite[Lemma 4.3]{CKS}; we omit the details. 
\end{proof} 

We are ready to prove Proposition \ref{prop:pos}. 

\begin{proof}[Proof of Proposition \ref{prop:pos}]
Let $\bar\e$ be as in Lemma \ref{lem:pos}, and suppose that $c\e \leq \bar \e$. Using $|\{ u \leq\frac{1}{2}\}\cap Q_1|\leq \e$, we may follow the proof of \cite[Proposition 4.2]{CKS} to find a constant $\rho$, depending only on $n$, $p$ and $q$, such that 
\begin{equation}\label{eq:u-Lpq-rl}
\mean{Q_{4\rho}} \left(\frac{(u^-)^q}{\rho^q} + |Du^+|^p\right) dx \leq c \e \rho^{q - p}\mean{Q_{4\rho}} (u^+)^p\,dx. 
\end{equation}
Therefore, defining $u_\rho : Q_4\to \R$ by 
$$
u_\rho (x) = \frac{u^+(\rho x)}{(4\rho)^{-\frac{n}{p}} \|u^+\|_{L^p(Q_{4\rho})}} - \frac{u^-(\rho x)}{4^{-\frac{n}{q}}\rho^{1-\frac{p}{q} - \frac{n}{q}} \| u^+ \|_{L^p(Q_{4\rho})}^{\frac{p}{q}}},
$$
we see that $u_\rho \in W^{1,p\wedge q}(Q_4)$ is a local $(1+ \e)$-minimizer of the functional $J$, such that
$$
\mean{Q_4} (u_\rho^+)^p\,dx = 1,\quad \mean{Q_4} ( (u_\rho^-)^q + |Du_\rho^+|^p) \,dx \leq c\e. 
$$
Since $c\e\leq \bar\e$, with $\bar\e$ as in Lemma \ref{lem:msr}, we obtain 
$$
u_\rho > 0 \quad\text{a.e.\ in }Q_1. 
$$
Rescaling back, we obtain that $u>0$ a.e.\ in $Q_{4\rho}$ as desired. 
\end{proof}


\section{H\"older regularity}\label{sec:Ca}

In this section, we study the universal H\"older regularity of local $(1+\e)$-minimizers for the functional $J_{p,q}$, and prove our first main result, Theorem \ref{thm:Ca}. Let us begin with a lemma that tells us how each phase of local minimizers for the functional $J_{p,q}$ should scale relatively to one another.

\begin{lem}\label{lem:comp}
Let $u\in W^{1,p\wedge q}(Q_1)$ be a local minimizer of the functional $J$, such that $\| u^+ \|_{L^p(Q_1)} = 1$ and $u(0) = 0$. If $\| u^+\|_{L^p(Q_{1/2})} \geq \beta$ for some $\beta>0$, then $\| u^- \|_{L^q(Q_1)} \geq c_\beta$, for some positive constant $c_\beta$ depending only on $n$, $p$, $q$ and $\beta$. 
\end{lem}

\begin{proof}
Let $\beta$ be any constant, with $0< \beta < 1$. Assume by way of contradiction that there exists a minimizer $u_j\in W^{1,p\wedge q}(Q_1)$ of the functional $J$, such that $\|u_j^+ \|_{L^p(Q_1)} = 1$, $\|u_j^+\|_{L^p(Q_{1/2})} \geq \beta$, $u_j(0) = 0$ but $\|u_j^-\|_{L^q(Q_1)} \leq \frac{1}{j}$. By \cite[Theorem 1.2]{CKS}, $u_j\in C^{0,\sigma}(Q_{1/2})$ and $\|u_j^+\|_{C^{0,\sigma}(Q_{1/2})} \leq c\| u_j^+\|_{L^p(Q_1)} \leq c$, and similarly, $\|u_j^-\|_{C^{0,\sigma}(Q_{1/2})} \leq \frac{c}{j}$, where both $c$ and $\sigma$ depend only on $n$, $p$ and $q$. This together with the Cacciopoli inequality (Lemma \ref{lem:caccio} with $\e= 0$) implies that $u_j^+ \to u_0$ weakly in $W^{1,p}(Q_{1/2})$ and uniformly in $Q_{1/2}$, while $u_j^- \to 0$ weakly in $W^{1,q}(Q_{1/2})$ and uniformly in $Q_{1/2}$, for some nonnegative function $u_0\in W^{1,p}(Q_{1/2})$. The uniform convergence along with $u_j(0) = 0$ implies that $u_0( 0) = 0$. In addition, passing to the limit in $\|u_j^+\|_{L^p(Q_{1/2})} \geq \beta$ ensures that $\|u_0 \|_{L^p(Q_{1/2})} \geq \beta$. However, the weak convergence of the gradient of $u_j$ implies that $u_0$ is also a minimizer of the functional $J$. As $u_0\geq 0$ in $Q_{1/2}$, $u_0$ is a $p$-harmonic function, but then it violates the minimizer principle, as $\|u_0 \|_{L^p(Q_{1/2})} \geq \beta>0$. 
\end{proof}

\begin{lem}\label{lem:comp-gr} 
Let $u\in W^{1,p\wedge q}(Q_1)$ be a local minimizer of the functional $J$, such that 
$$
\| u^+ \|_{L^p(Q_1)} \leq 1,\quad u(0) = 0,\quad \sup_{0<r<1} \frac{1}{r^{n+\sigma_- q}} \int_{Q_r} (u^-)^q\,dx \leq c_-,
$$
for some constants $c_-> 0$ and $\sigma_-\in(0,1]$. Then with $\sigma_+ = 1 - (1-\sigma_-) \frac{q}{p}$, 
$$
\sup_{0<r<1}  \frac{1}{r^{n + \sigma_+ p} } \int_{Q_r} (u^+)^p\,dx \leq c_+,
$$
where $c_+$ depends only on $n$, $p$, $q$, $\sigma_-$ and $c_-$. 
\end{lem}

\begin{proof}
Let $c_-$, $\sigma_-$ be given, and set $\sigma_+$ as in the statement. Suppose that the conclusion of this lemma is false. Then for each $j\in\N$, one can find a minimizer $u_j\in W^{1,p\wedge q}(Q_1)$ for the functional $J$, such that 
\begin{equation}\label{eq:uj-comp-gr}
\int_{Q_1} (u_j^+)^p \,dx  \leq 1,\quad u_j(0) =0,\quad\sup_{0<r<1} \frac{1}{r^{n+\sigma_- q}} \int_{Q_r} (u_j^-)^q\,dx \leq c_-,
\end{equation}
but 
\begin{equation}\label{eq:Sj-gr}
S_j := \sup_{\frac{r_j}{2} \leq r\leq 1} \frac{1}{r^{n + \sigma_+ p}} \int_{Q_r} (u_j^+)^p\,dx \to \infty
\end{equation}
where the supremum is achieved at $r=\frac{1}{2} r_j$; since $\| u_j^+\|_{L^p(Q_1)} \leq 1$, we must have $r_j\to 0$. Define 
$$
v_j(y) := \frac{u_j^+(r_j y)}{r_j^{-\frac{n}{p}} \| u_j^+ \|_{L^p(Q_{r_j})}} - \frac{u_j^-(r_j y)}{r_j^{1 - \frac{p}{q} - \frac{n}{q}}\| u_j^+ \|_{L^p (Q_{r_j})}^{\frac{p}{q}}}. 
$$
Then $v_j$ is a local minimizer of the functional $J$, such that by \eqref{eq:uj-comp-gr} and \eqref{eq:Sj-gr}, $\|v_j^+ \|_{L^p(Q_1)} =1$, $ \| v_j^+\|_{L^p(Q_{1/2})} \geq 2^{-\sigma_+}$ and $v_j (0) = 0$. Therefore, Lemma \ref{lem:comp} yields that $\|v_j^-\|_{L^q(Q_1)} \geq c_{\sigma_+}$. This implies that
\begin{equation}\label{eq:uj--gr}
\begin{aligned}
\frac{1}{r_j^n } \int_{Q_{r_j}} (u_j^-)^q\,dx& \geq c_{\sigma_+}^q r_j^{q - p} \int_{Q_{r_j}} (u_j^+)^p\,dx\geq \frac{S_j^{\frac{q}{p}}c_{\sigma_+}^q }{2^{\sigma_+ q}}  r_j^{q-p + \sigma_+ p}.
\end{aligned}
\end{equation}
Putting \eqref{eq:uj-comp-gr} and \eqref{eq:uj--gr} together, and recalling that $\sigma_+ = 1-(1-\sigma_-)\frac{q}{p}$, 
$$
c_{\sigma_-}^q \geq \frac{S_j^{\frac{q}{p}} c_{\sigma_+}^q}{2^{\sigma_+q}},
$$
a contradiction to the assumption that $S_j\to \infty$. 
\end{proof} 

Thanks to the above lemma, we can prove Theorem \ref{thm:Ca} for minimizers of the functional $J$.

\begin{proof}[Proof of Theorem \ref{thm:Ca} for minimizers]
It suffices to consider the case $p > q$, and 
$$
\int_{Q_1} ((u^+)^p + (u^+)^q)\,dx = 1. 
$$
By \cite[Theorem 1.1]{CKS}, we already know that $u\in C^{0,\sigma}(Q_1)$ and that $[u]_{C^{0,\sigma}(Q_1)} \leq c$, where both $c > 0$ and $\sigma\in(0,1)$ depend only on $n$, $p$ and $q$. Hence, if $u(z) = 0$ at some $z\in Q_{1/2}$, then 
$$
\sup_{0<r<\frac{1}{2}}\frac{1}{r^{n + \sigma q}} \int_{Q_r(z)} (u^-)^q\,dx \leq c,
$$
which along with Lemma \ref{lem:comp-gr} implies that 
$$
\sup_{0<r<\frac{1}{2}} \frac{1}{r^{n+ p - q + \sigma q}} \int_{Q_r(z)} (u^+)^p\,dx \leq c, 
$$
where the constant $c$ in both displays depends only on $n$, $p$ and $q$. Since $p > q$, $1- (1-\sigma)\frac{q}{ p} > \sigma > 0$. Now setting $\sigma_- = \sigma$ and $\sigma_+ = 1-(1-\sigma) \frac{q}{p}$, we immediately verify the relation required between $\sigma_+$ and $\sigma_-$. Since the above growth estimates hold uniformly around all $z\in\{u = 0\}\cap Q_1$, and since $\Delta_p u = 0$ in $\{u>0\}\cap Q_1$ and $\Delta_q u = 0$ in $\{ u < 0\}\cap Q_1$, one may arrive at the conclusion via some standard manipulation. We skip the detail. 
\end{proof}

Given a measurable function $u:\Omega\to\R$, define $D^+(u)$, $D^-(u)$ and $\Gamma(u)$ by the subset of $\Omega$ as follows:
$$
D^+(u) = \{ z\in \Omega: u > 0\text{ a.e.\ in some }Q_r(z)\subset \Omega\},\quad D^-(u) = D^+(-u),
$$
and 
$$
\Gamma(u) = \Omega\setminus (D^+(u)\cup D^-(u)). 
$$
By definition, both $D^+(u)$ and $D^-(u)$ are open and hence $\Gamma(u)$ is closed (relative to the topology of $\Omega$). Moreover, $z\in\Gamma(u)$ if and only if $|\{ u \geq 0\}\cap Q_r(z)| |\{ u \leq 0\}\cap Q_r(z)| > 0$ for any cube $Q_r(z)\subset \Omega$.

With Proposition \ref{prop:pos} at hand, we shall obtain, as a contraposition along with Lemma \ref{lem:Ca} below, that if a local $(1+\e)$-minimizer vanishes (in an appropriate Lebesgue sense) at certain point in the interior, then each phase exhibits certain universal H\"older growth. More exactly, we assert the following.

\begin{prop}\label{prop:Ca}
There exists a constant $\bar\sigma\in(0,1)$, depending only on $n$, $p$, and $q$, for which the following holds: for each $\sigma\in(0,\bar\sigma)$, one can find a constant $\e_\sigma \in (0,1)$, depending only on $n$, $p$, $q$, and $\sigma$, such that if $u\in W^{1,p\wedge q}(Q_1)$ is a local $(1+\e_\sigma)$-minimizer of the functional $J$ satisfying 
$$
\int_{Q_1} ((u^+)^p + (u^-)^q)\,dx \leq 1,\quad 0 \in \Gamma(u),
$$
then with $\sigma_+ = \sigma$ and $\sigma_- = 1-(1-\sigma)\frac{p}{q}$, one has
$$
\sup_{0<r<1}\left[ \frac{1}{r^{\sigma_+ p}} \mean{Q_r} (u^+)^p \,dx + \frac{1}{r^{ \sigma_- q}} \mean{Q_r} (u^-)^q\,dx \right] \leq c_\sigma, 
$$
where $c_\sigma$ depends only on $n$, $p$, $q$, and $\sigma$. 
\end{prop}

The following lemma will play a key role (together with Proposition \ref{prop:pos}). 

\begin{lem}\label{lem:Ca}
There exists a constant $\bar\sigma\in(0,1)$, depending only on $n$, $p$, and $q$, for which the following holds: for each $\sigma\in(0,\bar\sigma)$ and each $\tau \in (0,\frac{1}{2}]$, one can find $\e_{\sigma,\tau}\in(0,1)$, depending only on $n$, $p$, $q$, $\sigma$, and $\tau$, such that if $u\in W^{1,p\wedge q}(Q_1)$ is a local $(1+\e_{\sigma_+,\tau})$-minimizer of the functional $J$ satisfying, 
\begin{equation}\label{eq:Ca-asmp}
\int_{Q_1} ((u^+)^p + (u^-)^q)\,dx  = 1, \quad \frac{ |E^+(u,Q_r)|}{|Q_r|}\wedge\frac{|E^-(u,Q_r)|}{|Q_r|} \geq \tau,
\end{equation}
for some $r\in(0,1)$, where 
$$
\begin{aligned}
E^+(u,Q_r) &= \left\{ u \leq \frac{1}{2} r\Lambda(u,Q_r)^{\frac{1}{p}} \right\}, \\
E^-(u,Q_r) &= \left\{ u \geq - \frac{1}{2}r \Lambda(u,Q_r)^{\frac{1}{q}} \right\}, \\
\Lambda(u,Q_r) &=  \mean{Q_r} \bigg( \frac{(u^+)^p}{r^p} + \frac{(u^-)^q}{r^q} \bigg)\,dx, 
\end{aligned}
$$
then with $\sigma_+ = \sigma$ and $\sigma_- = 1- (1-\sigma)\frac{p}{q}$, one has
$$
\sup_{r \leq \rho\leq 1}\left[ \frac{1}{\rho^{ \sigma_+ p}}\mean{Q_\rho} (u^+)^p\, dx  + \frac{1}{\rho^{ \sigma_-q}} \mean{Q_\rho} (u^-)^q\,dx \right]\leq c_{\sigma,\tau},
$$
where $c_{\sigma,\tau}$ depends on the same parameters that determine $\e_{\sigma,\tau}$. 
\end{lem}

\begin{proof}
Let $\bar\sigma$ be determined later, and fix $\tau\in(0,\frac{1}{2}]$, $\sigma\in(0,\bar\sigma)$, and set $\sigma_\pm$ as in the stastement. Suppose by way of contradiction that for each $j\in\N$, we can find some positive constant $\e_j\to 0$, some local $(1+\e_j)$-minimizer $u_j \in W^{1,p\wedge q}(Q_1)$ of the functional $J$, and a radius $r_j \in(0,1)$ such that 
\begin{equation}\label{eq:uj-Lp-msr}
\int_{Q_1} ((u_j^+)^p + (u_j^-)^q)\,dx \leq 1,\quad  \frac{ |E^+(u_j, Q_{r_j})|}{|Q_{r_j}|}\wedge\frac{|E^-(u_j, Q_{r_j})|}{|Q_{r_j}|} \geq \tau,,
\end{equation}
but
\begin{equation}\label{eq:uj-Lp-Lq}
S_j = \sup_{r_j \leq r \leq 1} \left( \frac{1}{r^{ \sigma_+ p}}\mean{Q_r} (u_j^+)^p\,dx+ \frac{1}{r^{ \sigma_- q}} \mean{Q_r} (u_j^-)^q\,dx\right) \to \infty, 
\end{equation}
with the supremum achieved at level $r=r_j$. In order for \eqref{eq:uj-Lp-Lq} to be compatible with the first equality in \eqref{eq:uj-Lp-msr}, we must have $r_j\to 0$. 

Define $v_j  : Q_{r_j^{-1}} \to \R$ by 
$$
v_j(y) = \frac{u_j^+(r_j y)}{S_j^{\frac{1}{p}} r_j^{\sigma_+}} - \frac{u_j^- (r_j y)}{S_j^{\frac{1}{q}} r_j^{\sigma_-}}.
$$
By the way that it is rescaled, $v_j$ is a local $(1+\e_j)$-minimizer of the functional $J$ in $Q_{1/r_j}$. Moreover, by \eqref{eq:uj-Lp-msr} along with the relation $S_j r_j^{\sigma_+ p} = \Lambda_j r_j^p$ and $S_j r_j^{\sigma_- q} = \Lambda_j r_j^q$, where $\Lambda_j = \Lambda(u_j,Q_{r_j})$, 
\begin{equation}\label{eq:vj-msr}
\left|\left\{ |v_j|\leq \frac{1}{2} \right \}\cap Q_1 \right|\geq \tau,
\end{equation}
and by \eqref{eq:uj-Lp-Lq}, 
\begin{equation}\label{eq:vj-Lp-Lq}
\sup_{1\leq R\leq r_j^{-1}} \left[\frac{1}{R^{\sigma_+ p}}\mean{Q_R} (v_j^+)^p\,dy +  \frac{1}{R^{ \sigma_- q}} \mean{Q_R} (v_j^-)^q\,dy \right] = 1, 
\end{equation}
where the supremum is achieved at $R= 1$. 

Thanks to \eqref{eq:vj-msr} and \eqref{eq:vj-Lp-Lq}, one can argue analogously in the proof of Lemma \ref{lem:blow} to obtain a minimizer $v\in W_{loc}^{1,p\wedge q}(\R^n)$ of the functional $J$, with $v^+\in W_{loc}^{1,p}(\R^n)$ and $v^-\in W_{loc}^{1,q}(\R^n)$ such that
\begin{equation}\label{eq:v0-msr}
\left|\left\{ |v| \leq \frac{1}{2} \right \}\cap Q_1 \right|\geq \tau,
\end{equation}
and
\begin{equation}\label{eq:v0-Lp-Lq}
\sup_{R \geq 1}\left[ \frac{1}{R^{\sigma_+ p}} \mean{Q_R} (v^+)^p \, dy + \frac{1}{R^{ \sigma_- q}} \mean{Q_R} (v^-)^q \, dy \right] = 1,
\end{equation} 
where the supremum is achieved at $R = 1$. In particular, the latter observation indicates that $v$ is nontrivial. 

At this point, we choose $\bar\sigma$ as the positive exponent for Theorem \ref{thm:Ca} for minimizers; let us remind the readers that the statement for minimizers is proved right after the proof of Lemma \ref{lem:comp-gr}. Set $\bar\sigma_+ := \bar\sigma$ and $\bar\sigma_- := 1 - (1-\bar\sigma)\frac{q}{p}$. As $v$ is a local minimizer of the functional $J$ in $Q_{2R}$, $v^+\in C^{0,\bar \sigma_+}(Q_R)$, $v^-\in C^{0,\bar\sigma_-}(Q_R)$ and then by \eqref{eq:v0-Lp-Lq}, we derive that 
$$
[v^+]_{C^{0,\bar\sigma_+}(Q_R)}^p + [v^-]_{C^{0,\bar\sigma_-}(Q_R)}^q \leq \frac{c}{R^{(\bar\sigma_+ - \sigma_+)p}} + \frac{c}{R^{(\bar\sigma_- - \sigma_-)q}},
$$
for any $R>1$. As $\sigma_+ = \sigma < \bar\sigma_+$ and $\sigma_- = 1- (1-\sigma)\frac{q}{p} < \bar\sigma_-$, sending $R\to\infty$ implies that both $v^+$ and $v^-$ must be constant. Then by \eqref{eq:v0-msr}, $|v| \leq \frac{1}{2}$ everywhere in $Q_1$, whence $\int_{Q_1} ((v^+)^p + (v^+)^q)\,dx \leq 2^{-p} + 2^{-q} < 1$, a contradiction to the observation that the supremum in \eqref{eq:v0-Lp-Lq} is attained at $R=1$. 
\end{proof}

We are ready to prove Proposition \ref{prop:Ca}

\begin{proof}[Proof of Proposition \ref{prop:Ca}]
As $0\in\Gamma(u)$, there are three cases to consider: (i) $|\{ u >0\}\cap Q_\rho| |\{ u <0\}\cap Q_\rho| > 0$ for all $\rho\in(0,1)$, (ii) $u \geq 0$ a.e.\ in $Q_\rho$ for some small $\rho>0$, and (iii) $u\leq 0$ a.e.\ in $Q_\rho$ for some small $\rho> 0$. The last two cases are symmetric, and in those cases $u$ becomes a local $(1+\e)$-minimizer for the functional $J_{p,p}$, or $J_{q,q}$ depending on its sign. Thus, the growth estimate follows easily, once we establish the estimate for the first case. We leave out this part as an exercise for the reader. 

Henceforth, let us assume that the first case holds. Let $(\e,\tau,\mu)$ be the triple of constants from Proposition \ref{prop:pos} that  are determined solely by $n$, $p$ and $q$. Fix any $r\in(0,1)$. Since $|\{u>0\}\cap Q_{\mu r}|\cdot |\{ u< 0\}\cap Q_{\mu r}| > 0$, as a contraposition (applied to both $u$ and $-u$, after suitable rescaling), we obtain that 
\begin{equation}\label{eq:Ca-re}
\frac{ |E^+(u,Q_r)|}{|Q_r|}\wedge\frac{|E^-(u,Q_r)|}{|Q_r|} \geq \tau,
\end{equation}
with $E^+(u,Q_r)$ and $E^-(u,Q_r)$ defined as in Lemma \ref{lem:Ca}. As $\tau$ being a constant depending only on $n$, $p$ and $q$, the conclusion of this proposition follows immediately from Lemma \ref{lem:Ca}; this final step introduces another condition on the size of $\e$, which through the dependence of $\tau$ would be determined again solely by $n$, $p$, $q$, and $\sigma$. 
\end{proof}


\section{Almost Lipschitz regularity}\label{sec:C01-}

Here we prove almost Lipschitz regularity of almost minimizers to $J = J_{p,q}$, when $p$ and $q$ are close. 
Our proof is based on the compactness argument. The basic ingredient is the universal H\"older estimate for local minimizers of the functional $J_{p,q}$, see \cite[Theorem 1.2]{CKS}. 
Although it is not specified in the statement, one can observe from the higher integrability of each phase that the H\"older regularity is uniform when $p$ (or $q$) is close to $n$. We record this fact as a lemma below, as the proof of \cite[Theorem 1.2]{CKS} makes use of the local boundedness and the Harnack inequality for weak $p$-harmonic functions, and the constants involved in the latter assertions may vary as $p\to n$. 

\begin{lem}\label{lem:Ca-uni}
Let $u\in W^{1,p_+\wedge p_-}(Q_2)$ be a local minimizer of $J_{p_+,p_-}$. There exist $\bar\sigma \in (0,1)$, $c> 1$ and $\bar\delta>0$, all depending only on $n$, such that if $|n - p_\pm| \leq \bar\delta$, then 
$$
[u^\pm]_{C^{0,\bar\sigma}(Q_1)} \leq \bar c \| u^\pm \|_{L^{p_\pm}(Q_2)}.
$$
\end{lem}

\begin{proof}
Since $u$ is a local minimizer (instead of $(1+\e)$-minimizer) of $J_{p_+,p_-}$, $u^\pm$ is a weak $p_\pm$-subsolution in $Q_2$, according to \cite[Lemma 3.4]{CKS}. Hence, by \cite[Corollary 4.2]{GG}, there exist constants $\bar\delta > 0$ and $\bar\gamma \in (0,1)$, both depending only on $n$, such that if $|p_\pm - n| < \bar\delta$, then $u^\pm \in W^{1,p_\pm +  \bar\delta}(Q_1)\subset W^{1,n+ \bar\gamma \bar\delta}(Q_1)$. Now setting $\bar\sigma := 1- \frac{n}{n +\bar\gamma \bar\delta}$, it follows from the Sobolev embedding, the higher integrability and the Cacciopoli inequality for weak $p_\pm$-subsolutions that 
$$
\begin{aligned}
[u^\pm]_{C^{0,\bar\sigma}(Q_1)} &\leq c_1(n)\left[ \int_{Q_1} |Du^\pm|^{n + \bar\gamma \bar\delta}\,dx \right]^{\frac{1}{n+\bar\gamma \bar\delta}} \\
&  \leq c_1(n)c_2(n,p_\pm) \left[ \int_{Q_{3/2}} |Du^\pm|^{p_\pm}\,dx \right]^{\frac{1}{p_\pm}} \\
& \leq c_1(n)c_2(n,p_\pm) c_3(n,p_\pm) \left[ \int_{Q_2} (u^\pm)^{p_\pm}\,dx \right]^{\frac{1}{p_\pm}}. 
\end{aligned}
$$
Note that $c_2(n,p_\pm)$, and $c_3(n,p_\pm)$ are constants from the higher integrability and respectively the Cacciopoli inequality, and these are all uniformly bounded by a constant $c(n)$, as $p\to p_\pm$. Hence, our proof is finished. 
\end{proof}

Let us first verify the uniform growth of order $\sigma$ at free boundary points for minimizers. We prove it by compactness. 

\begin{lem}\label{lem:C01-}
Let $u\in W^{1,p\wedge q}(Q_1)$ be a local minimizer of $J_{p,q}$ such that 
\begin{equation}\label{eq:C01--asmp}
\int_{Q_1} ((u^+)^p + (u^-)^q)\,dx \leq 1,\quad u(0) = 0.
\end{equation}
Then for any $\sigma\in(0,1)$, there exists $\delta>0$, depending only on $n$, $p$, and $\sigma$, such that if $|p-q|<\delta$, then with $\sigma_+ = \sigma$ and $\sigma_- = 1-(1-\sigma)\frac{p}{q}$, 
$$
\frac{1}{r^{n+\sigma_+p}} \int_{Q_r} (u^+)^p \, dx + \frac{1}{r^{n+\sigma_- q}} \int_{Q_r} (u^-)^q \, dx \leq c,\quad\forall r\in(0,1),
$$
where $c>1$ depends only on $n$, $p$, and $\sigma$.
\end{lem}

\begin{proof}

Let $\sigma>0$ and $p\in(1,\infty)$ be given. Suppose that the conclusion of this lemma does not hold. Then for each $j=1,2,\cdots$, there must exist an exponent $q_j  > 1$ with $|q_j - p| \searrow 0$, a local minimizer $u_j \in W^{1,p\wedge q_j}(Q_1)$ of the functional $J_{p,q_j}$, and a scale $r_j \in (0,1)$, such that 
\begin{equation}\label{eq:uj-Linf-0}
\int_{Q_1} ((u_j^+)^p + (u_j^-)^{q_j}) \,dx \leq 1,\quad u_j(0) = 0,
\end{equation}
but with $\sigma_+ = \sigma$ and $\sigma_{j,-} = 1 - (1-\sigma)\frac{p}{q_j}\to \sigma$, 
\begin{equation}\label{eq:Sj-C01-}
S_j := \sup_{r_j \leq r\leq 1} \left[ \frac{1}{r^{\sigma_+ p}} \mean{Q_r}\,dx  +\frac{1}{r^{\sigma_{j,-} q_j}} \mean{Q_r} (u_j^-)^{q_j}\right] dx  \nearrow \infty.  
\end{equation}
To have the first inequality in \eqref{eq:uj-Linf-0} and \eqref{eq:Sj-C01-} to be compatible, we must have $r_j\searrow 0$ up to a subsequence. As in the proof of Lemma \ref{lem:Ca}, we consider the rescaling
$$
v_j (y) := \frac{u_j^+(r_j y)}{S_j^{\frac{1}{p}} r_j^{\sigma_+}} - \frac{u_j^-(r_j y)}{S_j^{\frac{1}{q_j}} r_j^{\sigma_{j,-}}}. 
$$
Then $v_j$ is a minimizer of $J_{p,q_j}$ in $Q_{1/r_j}$ and that 
\begin{equation}\label{eq:vj-C01-}
\sup_{1\leq R\leq\frac{1}{r_j}} \left[ \frac{1}{R^{\sigma_+p}}\mean{Q_R} (v_j^+)^p \,dx +\frac{1}{R^{\sigma_{j,-}q_j}} \mean{Q_R} (v_j^-)^{q_j} \,dx \right]= 1. 
\end{equation}
Then by \cite[Theorem 1.2]{CKS}, we have
\begin{equation}\label{eq:vj-C0a}
\sup_j \| v_j \|_{C^{0,\bar\sigma}(Q_R)} <  \infty, 
\end{equation}
where  both $c>1$ and $\bar\sigma\in(0,1)$ depend only on $n$ and $p$; see Lemma \ref{lem:Ca-uni} for the stability of $\bar\sigma$ and $c$ for the case $p = n$.  Moreover,  by \cite[Lemma 3.4]{CKS}, $v_j^+$ and $v_j^-$ are respectively, weak $p$-  and  $q_j$-subsolution, so the higher integrability \cite[Theorem 4.1]{GG} applies. Utilizing $|q_j - p|\searrow 0$, there exists $\eta  > 0$, depending only on $n$ and $p$, such that 
\begin{equation}\label{eq:vj-W1pe}
\sup_j \mean{Q_R} |Dv_j|^{p + \eta}\,dx < \infty.
\end{equation}
Also observe from \eqref{eq:uj-Linf-0} that
\begin{equation}\label{eq:vj-0}
v_j (0) = 0.
\end{equation}

By \eqref{eq:vj-C0a} and \eqref{eq:vj-W1pe}, we can extract a subsequence of $\{v_j\}_{j=1}^\infty$ along which $v_j \to v$ weakly in $W_{loc}^{1,p+\eta}(\R^n)$ and locally uniformly in $\R^n$, for some $v\in W_{loc}^{1,p+\eta}\cap C_{loc}^{0,\sigma}(\R^n)$. Let us continue to denote this subsequence by $\{v_j\}_{j=1}^\infty$. The uniform convergence along with \eqref{eq:vj-0} implies that 
\begin{equation}\label{eq:w-0}
v(0) = 0.
\end{equation}
We claim that $v$ is a (weak) $p$-harmonic function in $\R^n$. 

For any large $j$, we have $q_j \in (p - \eta,p+\eta)$. 
By \eqref{eq:vj-W1pe}, there exists a function $v \in W_{loc}^{1,p+\eta}(\R^n)$ such that $Dv_j\to Dv$ weakly in $L_{loc}^{p+\eta}(\R^n;\R^n)$ and $v_j \to v$ strongly in $L_{loc}^{p+\eta}(\R^n)$. Now let $R > 0$ be given and $w \in W_{loc}^{1,p+\eta}(\R^n)$ be given such $\supp(v-w)\subset B_R$. Then by the Fubini theorem, we can find $\delta_j\to 0$ and $\kappa \in (1,2)$ such that 
$$
\int_{\partial B_{\kappa R}} \bigg(|Dv_j|^{p+\eta} + |Dw|^{p+\eta} + \frac{|v_j - w|^{p+\eta}}{\delta_j^{p+\eta}} \bigg)\,d\sigma \leq c. 
$$
Therefore, we can find an auxiliary function $\vp_j\in W_{loc}^{1,p+\eta}(\R^n)$ (c.f. \cite[Lemma 1]{Lu}, which actually covers vectorial maps and yields stronger estimate than what we use here)  with a small exponent $\alpha \in(0,1)$, such that 
$$
\vp_j(x)  = \begin{cases}
    v_j (x) & \text{if }|x| \geq \kappa R \\
    w ((1-\delta_j^\alpha)^{-1}x) & \text{if } |x| \leq (1-\delta_j^\alpha) \kappa R,
\end{cases}
$$
and $\vp_j$ is defined in the ring $B_{\kappa R}\setminus B_{(1-\delta_j^\alpha)\kappa R}$ in a way that 
$$
\int_{B_{\kappa R}\setminus B_{(1-\delta_j^\alpha)\kappa R}} |D\vp_j|^{p+\eta}\,dx \leq c \delta_j^\alpha \to 0;
$$
here $c$ in the last displayed formula  may depend on $\kappa$ and $R$ but not on $j$. Then since $|p - q_j| < \frac{1}{2}\eta$ for all large $j$, denoting  $p_j^+ = p$ and $p_j^- = q_j$, it follows from the last displayed formula  and the H\"older inequality that 
$$
\int_{B_{\kappa R}\setminus B_{(1-\delta_j^\alpha)\kappa R}} |D\vp_j^\pm|^{p_j^\pm}\,dx \leq c\delta_j^{\frac{\alpha}{p+\eta}p_j^\pm} \to 0. 
$$
Therefore, we obtain 
$$
\begin{aligned}
J_{p,q_j}(\vp_j,B_{\kappa R}) & = (1-\delta_j^\alpha)^n J_{p,q_j}(w,B_{\kappa R}) + o(1) \leq \int_{B_{\kappa R}} |Dw|^p\,dx + o(1). 
\end{aligned}
$$ 
This, combined with the minimality of $v_j$ and the lower semicontinuity of the functional (noting that since $J_{p,q_j}(\phi) \leq \liminf_{k\to\infty} J_{p,q_j}(\phi_k)$,  for each $j$ fixed, whenever $\phi_k\to \phi$ weakly in $W^{1,p+\eta}$, and $J_{p,q_j}(\phi) \to J_{p,p}(\phi)$, we have $J_{p,p}(\phi) \leq \liminf_{j\to\infty} J_{p,q_j}(\phi_j)$ as $j\to\infty$ by the diagonal argument, whenever $\phi_j\to \phi$ weakly in $W^{1,p+\eta}$) we have 
$$
\begin{aligned}
\int_{B_{\kappa R}} |Dv|^p\,dx &\leq \liminf_{j\to\infty} J_{p,q_j}(v_j,B_{\kappa R})  \\
& \leq \liminf_{j\to \infty} J_{p,q_j} (\vp_j,B_{\kappa R}) \leq \int_{B_{\kappa R}} |Dw|^p\,dx. 
\end{aligned}
$$
Especially, as $\supp (v-w)\Subset B_R$ and $\kappa > 1$, we have $v = w$ in $B_{\kappa R}\setminus B_R$, so 
$$
\int_{B_R} |Dv|^p \,dx \leq \int_{B_R} |Dw|^p\,dx. 
$$
Since $R> 0$ was arbitrary, and also the competitor $w \in W_{loc}^{1,p}(\R^n)$, we conclude that $v$ is a local minimizer of the $p$-Dirichlet energy in $\R^n$, i.e., $v$ is a weak $p$-harmonic function in $\R^n$. This verifies our claim.

Now letting $k\to\infty$ in \eqref{eq:vj-C01-} and using $q_j\to p$, we obtain
\begin{equation}\label{eq:w-C01-}
\sup_{R\geq 1} \frac{1}{R^{\sigma p}} \mean{Q_R} |v|^p \,dx = 1.
\end{equation}
By the interior Lipschitz estimate for $p$-harmonic functions, 
\begin{equation}\label{eq:estimate}
[v]_{C^{0,1}(Q_R)} \leq \frac{c}{R^{1-\sigma}},
\end{equation}
for some $c$ independent of $R$. Taking $R\to\infty$ in \eqref{eq:estimate},
 we derive that $v$ is constant in $\R^n$, which together with \eqref{eq:w-0} implies $v\equiv 0$. This is yields a contradiction against \eqref{eq:w-C01-}, and the proof is finished.
\end{proof}

Next we extend the above lemma to local $(1+\e)$-minimizers. 

\begin{lem}\label{lem:C01-almost}
For any $\sigma\in(0,1)$, there exists $\e,\delta>0$, depending only on $n$, $p$, and $\sigma$, such that for any $q \in (1,\infty)$ with $|p-q| <\delta$, and any local  $(1+\e)$-minimizer $u\in W^{1,p\wedge q}(Q_1)$ satisfying \eqref{eq:C01--asmp}, one has, with $\sigma_+ = \sigma$ and $\sigma_- = 1- (1-\sigma)\frac{p}{q}$, that
$$
\frac{1}{r^{n+\sigma_+p}}\int_{Q_r} (u^+)^p \,dx + \frac{1}{r^{n+\sigma_-q}} \int_{Q_r} (u^-)^q\,dx \leq c,\quad\forall r\in(0,1),
$$
where $c>1$ depends only on $n$, $p$, and $\sigma$.
\end{lem}

\begin{proof}
As already observed in the proof of Proposition \ref{prop:Ca}, the assumption $u(0) = 0$ implies \eqref{eq:Ca-re} for every $r\in(0,1)$. Hence, the assumption \eqref{eq:C01--asmp} implies \eqref{eq:Ca-asmp}. The rest of the proof is the same with that of Lemma \ref{lem:Ca}. More exactly, given $\sigma\in(0,1)$ and $p>1$, we first choose $\delta>0$ sufficiently small such that Lemma \ref{lem:C01-} holds with $\frac{1+\sigma}{2}$ in place of $\sigma$, for all local minimizers for functional $J_{p,q}$ for any $q\in(1,\infty)$ with $|p-q| < \delta$. Then we can take $\e>0$ small enough such that Lemma \ref{lem:Ca} holds with $\tau$ as in \eqref{eq:Ca-re}, $\sigma_+ = \sigma$ and $\sigma_- = 1- (1-\sigma)\frac{p}{q}$, $\bar\sigma_+ = \frac{1+\sigma}{2} > \sigma = \sigma_-$, and $\bar\sigma_- = 1 - (\frac{1-\sigma}{2})\frac{p}{q} > 1- (1-\sigma)\frac{p}{q} = \sigma_-$. We skip the details. 
\end{proof}

We are ready to prove the almost Lipschitz regularity for almost minimizers, when $|p-q|\ll 1$. 
\begin{proof}[Proof of Theorem \ref{thm:C01-}]
With the same (and simpler) compactness argument, we can also prove that local $(1+\e)$-minimizers for $J_{p,p} (w) \equiv \int |Dw|^p\,dx$ is of class $C^{0,\sigma}$, for any $\sigma\in(0,1)$ and every $\e \in (0, \e_\sigma)$, since $p$-harmonic functions are of class $C^{1,\alpha} \subset C^{0,1}$. Moreover, we can obtain a uniform $C^{0,\sigma}$-estimates, with this compactness argument, and the smallness constant $\e_\sigma$ depends only on $n$, $p$, and $\sigma$. Thus, the passage from Lemma \ref{lem:C01-almost} to Theorem \ref{thm:C01-} is standard. We shall not present the obvious details here. 
\end{proof}


\section*{Declarations}

\noindent {\bf  Data availability statement:} All data needed are contained in the manuscript.

\medskip
\noindent {\bf  Funding and/or Conflicts of interests/Competing interests:} The authors declare that there are no financial, competing or conflict of interests.


\end{document}